\documentclass[11pt]{article}

\usepackage[top=1.9in,bottom=1.9in,left=2in,right=2in]{geometry}
\usepackage{amsfonts} 
\usepackage{amsmath}
\usepackage{amssymb}
\usepackage{setspace}
\usepackage{subfigure}
\usepackage{url}
\usepackage{booktabs}
\usepackage{ifthen}
\usepackage{tikz}
\usetikzlibrary{decorations.pathreplacing}
\usetikzlibrary{calc}
\usepackage{enumerate}

\newcommand{\hs}{\hspace}

\newcommand{\qed}{\hs*{\fill} $\Box$}

\newtheorem{thm}{Theorem}[section]
\newtheorem{lem}[thm]{Lemma}

\newtheorem{ex}[thm]{Example}
\newtheorem{conj}[thm]{Conjecture}

\numberwithin{equation}{section}

\begin{document}

\title{On super edge-graceful trees of diameter four}

\author{
Elliot Krop\thanks{Department of Mathematics, Clayton State University, (\texttt{ElliotKrop@clayton.edu})}
\and
Fedelis Mutiso\thanks{Department of Mathematics, Clayton State University, (\texttt{fmutiso@student.clayton.edu})}
\and 
Christopher Raridan\thanks{Department of Mathematics, Clayton State University, (\texttt{ChristopherRaridan@clayton.edu})}}
   
\date{\today}

\maketitle

\begin{abstract}
In~\cite{CLGS}, \mbox{P.T. Chung}, \mbox{S.-M. Lee}, \mbox{W.Y. Gao}, and \mbox{K. Schaffer} posed the following problem: \textit{Characterize trees of diameter $4$ which are super edge-graceful}. In this paper, we provide super edge-graceful labelings for all caterpillars and even size lobsters of diameter $4$ which permit such labelings. We also provide super edge-graceful labelings for several families of odd size lobsters of diameter $4$.
\\[\baselineskip] 
\\[\baselineskip]
\end{abstract}

\section{Introduction}  

Let $G(V,E)$ be a finite, simple, undirected graph with vertex set $V$ and edge set $E$ such that the order is $|V|=p$ and the size is $|E|=q$. Edge $uv$ is incident with vertices $u$ and $v$. Rosa introduced the $\beta$-valuation of a graph in~\cite{ARosa}, which was later popularized as a graceful labeling by Golomb~\cite{GWGolomb}. A graph is \textit{graceful} if there exists an injection $f : V \rightarrow \{ 0, 1, \ldots q \}$ such that the induced edge labeling $f^* : E \rightarrow \{1, \ldots, q \}$ defined by $f^*(uv) = |f(u)-f(v)|$ is bijective. Lo~\cite{SPLo} defined a graph as \textit{edge-graceful} if the edges can be labeled from $\{ 1, \ldots, q \}$ and the resulting vertex sums are distinct $(\mbox{mod}~p)$. In~\cite{MS}, Mitchem and Simoson defined a graph to be \textit{super edge-graceful} if there exists a bijection $f : E \rightarrow \{ 0, \pm1, \ldots, \pm\frac{q-1}{2} \}$ when $q$ is odd and $f : E \rightarrow \{ \pm1, \ldots, \pm\frac{q}{2} \}$ when $q$ is even such that the induced vertex labeling $f^+$ for a vertex $v$, given by $f^+(v) = \sum_{uv \in E} f(uv)$, is a bijection $f^+ : V \rightarrow \{ 0, \pm1, \ldots, \pm\frac{p-1}{2} \}$ when $p$ is odd and $f^+ : V \rightarrow \{ \pm1, \ldots, \pm\frac{p}{2} \}$ when $p$ is even. Numerous papers have been written on the various graph labelings; see Gallian's dynamic survey~\cite{JAGallian} for more information. 

In 2006, Chung, Lee, Gao, and Schaffer posed the following problem in~\cite{CLGS}: Characterize trees of diameter $4$ which are super edge-graceful. In the following year, Lee and Ho proved all trees of odd order with three even vertices are super edge-graceful~\cite{LH}, which provided a partial solution to this problem. In this paper, we provide super edge-graceful labelings for all caterpillars and even size lobsters of diameter $4$ which permit such labelings. We also provide super edge-graceful labelings for several families of odd size lobsters of diameter $4$.  

Let $\{ a_i \}_{i=1}^n$ be a sequence of nonnegative integers, where at least two entries in the sequence are positive, that has been arranged so that all of the entries that are $0$ are at the beginning of the sequence, followed by all of the positive even entries (in nondecreasing order), followed by all of the positive odd entries (in nondecreasing order). That is, let $n=j+k+l$ and for $1 \leq i \leq j$, we have $a_i=0$; for $j+1 \leq i \leq j+k$, we have even $a_i > 0$ with $a_r \leq a_s$ for $j+1 \leq r<s \leq j+k$; and for $j+k+1 \leq i \leq n$, we have odd $a_i > 0$ with $a_r \leq a_s$ for $j+k+1 \leq r<s \leq n$. Thus, $j$ is the cardinality of $\{ a_i : a_i=0 \}$, $k$ is the cardinality of $\{ a_i : a_i>0,~\mbox{even} \}$, and $l$ is the cardinality of $\{ a_i : a_i>0,~\mbox{odd} \}$, where $k+l \geq 2$.

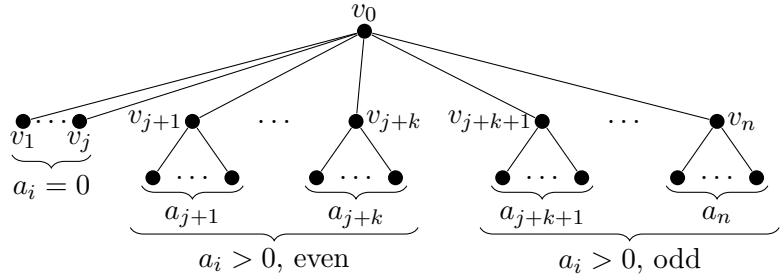
\begin{figure}[ht]
\centering
\begin{tikzpicture}[scale=1.5]
\tikzstyle{vertex1}=[circle,fill=black,inner sep=2pt]
\tikzstyle{vertex2}=[circle,fill=black,inner sep=0.5pt]

\node[vertex1] (c-0) at (3.025,1.3) {};
\coordinate[label=above:$v_0$] (C-0) at (3.025,1.3);

\foreach \cent/\lab/\xcoord in {1/1/0, j/j/0.5}
{
  \node[vertex1] (c-\cent) at (\xcoord,0.5) {};
  \coordinate[label=below:{$v_{\lab}$}] (C-\cent) at (\xcoord,0.5);
}

\draw[decorate,decoration={brace,amplitude=5pt},below=7pt] 
	(0.6,0.45) -- (-0.1,0.45) node[black,midway,below=4pt]{$a_i=0$};

\foreach \cent/\lab/\xcoord in {j+1/j+1/1.5, j+k+1/j+k+1/4.6}
{
  \node[vertex1] (c-\cent) at (\xcoord,0.5) {};
  \coordinate[label=left:{$v_{\lab}$}] (C-\cent) at (\xcoord,0.5);
}

\foreach \cent/\lab/\xcoord in {j+k/j+k/2.95, n/n/6.15}
{
  \node[vertex1] (c-\cent) at (\xcoord,0.5) {};
  \coordinate[label=right:{$v_{\lab}$}] (C-\cent) at (\xcoord,0.5);
}

\foreach \cent in {1,j,j+1,j+k,j+k+1,n}
 \draw (c-0) -- (c-\cent);
 
\foreach \xcoord in {0.25,2.225,5.325}
 	\node (dots) at (\xcoord,0.5) {$\ldots$}; 
 
\foreach \pen/\laba/\labb/\xcoord in {j+1-1/j+1/1/1.15, j+1-b/j+1/a_{j+1}/1.85, j+k-1/j+k/1/2.6, j+k-b/j+k/a_{j+k}/3.3, j+k+1-1/j+k+1/1/4.25, j+k+1-b/j+k+1/a_{j+k+1}/4.95, n-1/n/1/5.8, n-b/n/a_n/6.5}
{
  \node[vertex1] (x-\pen) at (\xcoord,0) {};
} 

\foreach \cent/\pen in {j+1/j+1-1, j+1/j+1-b, j+k/j+k-1, j+k/j+k-b, j+k+1/j+k+1-1, j+k+1/j+k+1-b, n/n-1, n/n-b}
 \draw (c-\cent) -- (x-\pen); 
 
\foreach \xcoord in {1.5,2.95,4.6,6.15}
 	\node (dots) at (\xcoord,0) {$\ldots$}; 
 	
\draw[decorate,decoration={brace,amplitude=5pt},below=2pt] 
	(1.95,0) -- (1.05,0) node[black,midway,below=4pt]{$a_{j+1}$};
	
\draw[decorate,decoration={brace,amplitude=5pt},below=2pt] 
	(3.4,0) -- (2.5,0) node[black,midway,below=4pt]{$a_{j+k}$};
	
\draw[decorate,decoration={brace,amplitude=5pt},below=2pt] 
	(5.05,0) -- (4.15,0) node[black,midway,below=4pt]{$a_{j+k+1}$};
	
\draw[decorate,decoration={brace,amplitude=5pt},below=2pt] 
	(6.6,0) -- (5.7,0) node[black,midway,below=4pt]{$a_{n}$};
 	
\draw[decorate,decoration={brace,amplitude=5pt},below=4pt] 
	(3.5,-0.3) -- (0.95,-0.3) node[black,midway,below=4pt]{$a_i>0$, even};
	
\draw[decorate,decoration={brace,amplitude=5pt},below=4pt] 
	(6.7,-0.3) -- (4.05,-0.3) node[black,midway,below=4pt]{$a_i>0$, odd};
\end{tikzpicture}
\caption{A general rooted tree of height $2$, $RT(a_1, \ldots, a_n)$.}
\label{fig:rooted-tree}
\end{figure}

Define a rooted tree of height $2$ (diameter $4$) as follows: let $v_0$ be the root and let the children of $v_0$ be $v_1, \ldots, v_n$ such that $v_i$ has $a_i$ children, where $a_i$ is defined by the above sequence. Then for $1 \leq i \leq j$, vertex $v_i$ has no children and for $j+1 \leq i \leq n$, the children of $v_i$ are labeled $v_{i,m}$, where $1 \leq m \leq a_i$. The edge $v_0v_r$ will be expressed as $e_{0,r}$ and the edge $v_rv_{r,s}$ will be expressed as $e_{r,s}$. Such a tree will be denoted $RT(a_1, a_2, \ldots, a_n)$, an example of which can be seen in Figure~\ref{fig:rooted-tree}. 
In $RT(a_1, a_2, \ldots, a_n)$, the number of edges is $q = n + \sum_{i=j+1}^n a_i$ and the number of vertices is $p = q+1$. For nonnegative integers $m$ and $n$, we use the abbreviated form $m^n$ to indicate that $m$ is written in a list $n$ times; that is, $RT(0^4,2^3)$ is $RT(0,0,0,0,2,2,2)$. If $n=0$, then $m$ is not written to the list.

\section{Caterpillars of diameter four}

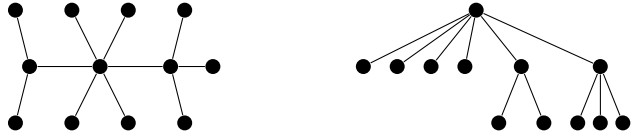
\begin{figure}[ht]
\centering
\begin{tikzpicture}[scale=1.5]
\tikzstyle{vertex1}=[circle,fill=black,inner sep=2pt]

\node[vertex1] (C-5) at (0.875,0.5) {};
\node[vertex1] (C-0) at (1.5,0.5) {};
\node[vertex1] (C-6) at (2.125,0.5) {};

\node[vertex1] (C-1) at (1.75,1) {};
\node[vertex1] (C-2) at (1.25,1) {};
\node[vertex1] (C-3) at (1.25,0) {};
\node[vertex1] (C-4) at (1.75,0) {};

\node[vertex1] (X-51) at (0.75,0) {}; 
\node[vertex1] (X-52) at (0.75,1) {};

\node[vertex1] (X-61) at (2.25,0) {};
\node[vertex1] (X-62) at (2.25,1) {};
\node[vertex1] (X-63) at (2.5,0.5) {};

\foreach \x in {1,2,...,6}
	\draw (C-0) -- (C-\x); 

\foreach \x in {51,52}
	\draw (C-5) -- (X-\x);

\foreach \x in {61,62,63}
	\draw (C-6) -- (X-\x);
	
\node[vertex1,xshift=5cm] (CC-5) at (1.9,0.5) {};
\node[vertex1,xshift=5cm] (CC-0) at (1.5,1) {};
\node[vertex1,xshift=5cm] (CC-6) at (2.6,0.5) {};

\node[vertex1,xshift=5cm] (CC-1) at (0.5,0.5) {};
\node[vertex1,xshift=5cm] (CC-2) at (0.8,0.5) {};
\node[vertex1,xshift=5cm] (CC-3) at (1.1,0.5) {};
\node[vertex1,xshift=5cm] (CC-4) at (1.4,0.5) {};

\node[vertex1,xshift=5cm] (XX-51) at (1.7,0) {}; 
\node[vertex1,xshift=5cm] (XX-52) at (2.1,0) {};

\node[vertex1,xshift=5cm] (XX-61) at (2.4,0) {};
\node[vertex1,xshift=5cm] (XX-62) at (2.6,0) {};
\node[vertex1,xshift=5cm] (XX-63) at (2.8,0) {};

\foreach \x in {1,2,...,6}
 	\draw (CC-0) -- (CC-\x); 

\foreach \x in {51,52}
	\draw (CC-5) -- (XX-\x);

\foreach \x in {61,62,63}
	\draw (CC-6) -- (XX-\x);	
	
\end{tikzpicture}
\caption{The caterpillar $C(2,4,3)$ is the rooted tree $RT(0^4,2,3)$.}
\label{fig:gen-cat}
\end{figure}

A caterpillar is a tree with the property that the removal of its endpoints leaves a path~\cite{HS}. All caterpillars can be represented as rooted trees and any caterpillar of diameter $4$ will be represented as a rooted tree of height $2$. For example, the caterpillar $C(2,4,3)$ is the rooted tree $RT(0^4,2,3)$ seen in Figure~\ref{fig:gen-cat}. 

Rooted trees that are caterpillars of diameter $4$, that is, \linebreak $RT(0^j, a_{j+1}, a_{j+2})$ for integers $j \geq 0$ and $a_{j+1},a_{j+2} \geq 1$, have $q = j + 2 + a_{j+1} + a_{j+2}$ edges and $p = q+1$ vertices. In the following two subsections, we examine $RT(0^j, a_{j+1}, a_{j+2})$ based on the parity of $q$, the size of the graph.

\subsection{Even size caterpillars of diameter four}

When $j$ is even and $a_{j+1},a_{j+2}$ are of the same parity, or when $j$ is odd and $a_{j+1},a_{j+2}$ are of opposite parity, $q$ is even and $p$ is odd. 

\begin{lem}
\label{lem:e,sameparity}
For all even integers $j \geq 0$ and integers $a_{j+2} \geq a_{j+1} \geq 1$ of the same parity, the rooted tree $RT(0^j, a_{j+1}, a_{j+2})$ is super edge-graceful.
\end{lem}

\begin{proof}
Let $j=2r$ for some integer $r \geq 0$. 

\textit{Case 1.} Suppose $a_{2r+1}=2s$ and $a_{2r+2}=2t$ for some integers $t \geq s \geq 1$. Then $G=RT(0^{2r}, 2s, 2t)$ has $2(r+s+t+1)$ edges and $2(r+s+t+1)+1$ vertices. An edge-labeling bijection $f$ maps $E(G)$ to $\{ \pm1, \ldots, \pm(r+s+t+1) \}$ as follows: for $1 \leq i \leq r+1$, $f(e_{0,2i-1}) = i$ and $f(e_{0,2i}) = -i$; for $1 \leq i \leq s$, $f(e_{2r+1,2i-1}) = r+1+i$ and $f(e_{2r+1,2i}) = -(r+1+i)$; and for $1 \leq i \leq t$, $f(e_{2r+2,2i-1}) = r+s+1+i$ and $f(e_{2r+2,2i}) = -(r+s+1+i)$. This edge labeling induces a vertex labeling $f^+$ from $V(G)$ to $\{ 0, \pm1, \ldots, \pm(r+s+t+1) \}$ as follows: $f^+(v_0)=0$; for $1 \leq i \leq 2r+2$, $f^+(v_i) = f(e_{0,i})$; for $1 \leq i \leq 2s$, $f^+(v_{2r+1,i}) = f(e_{2r+1,i})$; and for $1 \leq i \leq 2t$, $f^+(v_{2r+2,i}) = f(e_{2r+2,i})$. 

\textit{Case 2.} Suppose $a_{2r+1}=2s-1$ and $a_{2r+2}=2t-1$ for some integers $t \geq s \geq 1$. Then $G=RT(0^{2r}, 2s-1, 2t-1)$ has $2(r+s+t)$ edges and $2(r+s+t)+1$ vertices. An edge-labeling bijection $f$ maps $E(G)$ to $\{ \pm1, \ldots, \pm(r+s+t) \}$ as follows: $f(e_{0,2r+1}) = 1$, $f(e_{0,2r+2}) = -1$, $f(e_{2r+1,1}) = -2$, and $f(e_{2r+2,1}) = 2$; for $1 \leq i \leq r$, $f(e_{0,2i-1}) = 2+i$ and $f(e_{0,2i}) = -(2+i)$; for $1 \leq i \leq s-1$, $f(e_{2r+1,2i}) = r+2+i$ and $f(e_{2r+1,2i+1}) = -(r+2+i)$; and for $1 \leq i \leq t-1$, $f(e_{2r+2,2i}) = r+s+1+i$ and $f(e_{2r+2,2i+1}) = -(r+s+1+i)$. This edge labeling induces a vertex labeling $f^+$ from $V(G)$ to $\{ 0, \pm1, \ldots, \pm(r+s+t+1) \}$ as follows: $f^+(v_0)=0$; for $1 \leq i \leq 2r$, $f^+(v_i) = f(e_{0,i})$; $f^+(v_{2r+1}) = -1$ and $f^+(v_{2r+2}) = 1$; for $1 \leq i \leq 2s-1$, $f^+(v_{2r+1,i}) = f(e_{2r+1,i})$; and for $1 \leq i \leq 2t-1$, $f^+(v_{2r+2,i}) = f(e_{2r+2,i})$. \qed
\end{proof}

\begin{lem}
\label{lem:o,oppparity}
For all odd integers $j,a_{j+2} \geq 1$ and even integers \linebreak $a_{j+1} \geq 2$, the rooted tree $RT(0^j, a_{j+1}, a_{j+2})$ is super edge-graceful.
\end{lem}

\begin{proof}
Let $j=2r-1$, $a_{2r}=2s$, and $a_{2r+1}=2t-1$ for some integer $r,s,t \geq 1$. Then $G=RT(0^{2r-1}, 2s, 2t-1)$ has $2(r+s+t)$ edges and $2(r+s+t)+1$ vertices. An edge-labeling bijection $f$ maps $E(G)$ to $\{ \pm1, \ldots, \pm(r+s+t) \}$ as follows: $f(e_{0,2r+1}) = 1$ and $f(e_{2r+1,1}) = -1$; for $1 \leq i \leq r$, $f(e_{0,2i-1}) = 1+i$ and $f(e_{0,2i}) = -(1+i)$; for $1 \leq i \leq s$, $f(e_{2r,2i-1}) = r+1+i$ and $f(e_{2r,2i}) = -(r+1+i)$; and for $1 \leq i \leq t-1$, $f(e_{2r+1,2i}) = r+s+1+i$ and $f(e_{2r+1,2i+1}) = -(r+s+1+i)$. This edge labeling induces a vertex labeling $f^+$ from $V(G)$ to $\{ 0, \pm1, \ldots, \pm(r+s+t) \}$ as follows: $f^+(v_0)=1$; for $1 \leq i \leq 2r$, $f^+(v_i) = f(e_{0,i})$; $f^+(v_{2r+1})=0$; for $1 \leq i \leq 2s$, $f^+(v_{2r,i}) = f(e_{2r,i})$; and for $1 \leq i \leq 2t-1$, $f^+(v_{2r+1,i}) = f(e_{2r+1,i})$. \qed
\end{proof}
\\[\baselineskip]
\indent Lemma~\ref{lem:e,sameparity} and Lemma~\ref{lem:o,oppparity} imply the following theorem:

\begin{thm}
\label{thm:evencaterpillars}
All caterpillars of diameter $4$ with even size are super edge-graceful.
\end{thm}

\begin{ex}
\label{ex:even-size-cat}
Figure~\ref{fig:RT(0^4,2,6)-RT(0^3,2,5)} shows super edge-graceful labelings of two even size caterpillars of diameter $4$, $RT(0^4, 2, 6)$ and $RT(0^3,2,5)$. Edges incident with pendant vertices have the same label as those vertices.
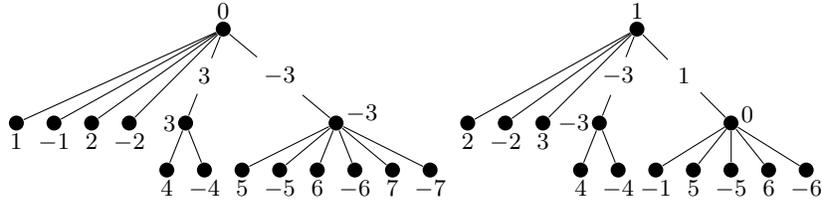
\begin{figure}[ht]
\centering
\begin{tikzpicture}[scale=2.5]
\tikzstyle{vertex}=[circle,fill=black,inner sep=2pt]


\node[vertex] (c-0) at (1.2,1) {};

\coordinate[label=above:{\footnotesize $0$}] (C-0) at (1.2,1);

\foreach \cent/\lab/\xcoord in {1/1/0.1, 2/-1/0.3, 3/2/0.5, 4/-2/0.7}
{
  \node[vertex] (c-\cent) at (\xcoord,0.5) {};
  \coordinate[label=below:{\footnotesize $\lab$}] (C-\cent) at (\xcoord,0.5);
}

\node (e-05) at (1.1,0.75) {\footnotesize $3$};
\node (e-06) at (1.5,0.75) {\footnotesize $-3$};

\node[vertex] (c-5) at (1,0.5) {};
\node[vertex] (c-6) at (1.8,0.5) {};

\coordinate[label=left:{\footnotesize $3$}] (C-5) at (1,0.5);
\coordinate[label=above right:{\footnotesize $-3$}] (C-6) at (1.8,0.45);

\foreach \pen/\lab/\xcoord in {51/4/0.9, 52/-4/1.1, 61/5/1.3, 62/-5/1.5, 63/6/1.7, 64/-6/1.9, 65/7/2.1, 66/-7/2.3}
{
  \node[vertex] (x-\pen) at (\xcoord,0.25) {};
  \coordinate[label=below:{\footnotesize $\lab$}] (X-\pen) at (\xcoord,0.25);
}

\draw 
(c-0) -- (e-05) -- (c-5)
(c-0) -- (e-06) -- (c-6);
\foreach \lab in {1,2,3,4}
 \draw (c-0) -- (c-\lab);
\foreach \lab in {51,52}
 \draw (c-5) -- (x-\lab);
\foreach \lab in {61, 62, 63, 64, 65, 66}
 \draw (c-6) -- (x-\lab);
 
 
\node[vertex,xshift=5.5cm] (c-0) at (1.2,1) {};

\coordinate[label=above:{\footnotesize $1$},xshift=5.5cm] (C-0) at (1.2,1);

\foreach \cent/\xcoord in {1/0.3, 2/0.5, 3/0.7, 4/1, 5/1.7}
 \node[vertex,xshift=5.5cm] (c-\cent) at (\xcoord,0.5) {};

\foreach \cent/\lab/\xcoord in {1/2/0.3, 2/-2/0.5, 3/3/0.7}
 \coordinate[label=below:{\footnotesize $\lab$},xshift=5.5cm] (C-\cent) at (\xcoord,0.5);
 
\coordinate[label=left:{\footnotesize $-3$},xshift=5.5cm] (C-4) at (1,0.5);
\coordinate[label=above right:{\footnotesize $0$},xshift=5.5cm] (C-5) at (1.7,0.45);

\node[xshift=5.5cm] (e-04) at (1.1,0.75) {\footnotesize $-3$};
\node[xshift=5.5cm] (e-05) at (1.45,0.75) {\footnotesize $1$};

\foreach \pen/\lab/\xcoord in {41/4/0.9, 42/-4/1.1, 51/-1/1.3, 52/5/1.5, 53/-5/1.7, 54/6/1.9, 55/-6/2.1}
{
  \node[vertex,xshift=5.5cm] (x-\pen) at (\xcoord,0.25) {};
  \coordinate[label=below:{\footnotesize $\lab$},xshift=5.5cm] (X-\pen) at (\xcoord,0.25);
}
 
\foreach \cent in {1,2,3}
 \draw (c-0) -- (c-\cent);
\draw
(c-0) -- (e-04) -- (c-4)
(c-0) -- (e-05) -- (c-5);
\foreach \lab in {41, 42}
 \draw (c-4) -- (x-\lab);
\foreach \lab in {51, 52, 53, 54, 55}
 \draw (c-5) -- (x-\lab); 

\end{tikzpicture}
\caption{$RT(0^4, 2, 6)$ and $RT(0^3,2,5)$ are super edge-graceful.}
\label{fig:RT(0^4,2,6)-RT(0^3,2,5)}
\end{figure}
\end{ex}

\subsection{Odd size caterpillars of diameter four}

When $j$ is even and $a_{j+1},a_{j+2}$ are of opposite parity, or when $j$ is odd and $a_{j+1},a_{j+2}$ are of the same parity, $q$ is odd and $p$ is even. If $G=RT(0^j, a_{j+1}, a_{j+2})$ is to be super edge-graceful, then one of the edges in $E(G)$ must be labeled by $0$. Since the edges in $E(G) - \{ e_{0,j+1}, e_{0,j+2} \}$ are incident to pendant vertices, none of these edges can be labeled by $0$, which means that either $f(e_{0,j+1})=0$ or $f(e_{0,j+2})=0$.

\begin{lem}
\label{lem:e,oppparity}
For all even integers $j \geq 0$, even integers $a_{j+1} \geq 2$, and odd integers $a_{j+2} \geq 1$, the rooted tree $RT(0^j, a_{j+1}, a_{j+2})$ is super edge-graceful.
\end{lem}

\begin{proof}
Let $j=2r$, $a_{2r+1}=2s$, and $a_{2r+2}=2t-1$ for some integers $r \geq 0$, $s,t \geq 1$. Then $G=RT(0^{2r}, 2s, 2t-1)$ has $2(r+s+t)+1$ edges and $2(r+s+t+1)$ vertices. An edge-labeling bijection $f$ maps $E(G)$ to $\{ 0, \pm1, \ldots, \pm(r+s+t) \}$ as follows: $f(e_{0,2r+1}) = 0$, $f(e_{0,2r+2}) = 1$, $f(e_{2r+1,1}) = -1$, $f(e_{2r+1,2}) = -(r+s+t)$, and $f(e_{2r+2,1}) = r+s+t$; for $1 \leq i \leq r$, $f(e_{0,2i-1}) = 1+i$ and $f(e_{0,2i}) = -(1+i)$; for $2 \leq i \leq s$, $f(e_{2r+1,2i-1}) = r+i$ and $f(e_{2r+1,2i}) = -(r+i)$; and for $1 \leq i \leq t-1$, $f(e_{2r+2,2i}) = r+s+i$ and $f(e_{2r+2,2i+1}) = -(r+s+i)$. This edge labeling induces a vertex labeling $f^+$ from $V(G)$ to $\{ \pm1, \ldots, \pm(r+s+t+1) \}$ as follows: $f^+(v_0)=1$; for $1 \leq i \leq 2r$, $f^+(v_i) = f(e_{0,i})$; $f^+(v_{2r+1}) = -(r+s+t+1)$ and $f^+(v_{2r+2}) = r+s+t+1$; for $1 \leq i \leq 2s$, $f^+(v_{2r+1,i}) = f(e_{2r+1,i})$; and for $1 \leq i \leq 2t-1$, $f^+(v_{2r+2,i}) = f(e_{2r+2,i})$. \qed
\end{proof}

\begin{lem}
\label{lem:o,e,e}
For all odd integers $j \geq 1$ and even integers $a_{j+2} \geq a_{j+1} \geq 2$, the rooted tree $RT(0^j, a_{j+1}, a_{j+2})$ is super edge-graceful.
\end{lem}

\begin{proof}
Let $l=2r-1$, $a_{2r}=2s$, and $a_{2r+1}=2t$ for some integers $r,s,t \geq 1$. Then $G=RT(0^{2r-1},2s,2t)$ has $2(r+s+t)+1$ edges and $2(r+s+t+1)$ vertices. An edge-labeling bijection $f$ maps $E(G)$ to $\{ 0, \pm1, \ldots, \pm(r+s+t) \}$ as follows: $f(e_{0,1})=1$, $f(e_{0,2r})=0$, $f(e_{0,2r+1}) = r+s+t$, $f(e_{2r,1})=-1$, and $f(e_{2r,2}) = -(r+s+t)$; for $1 \leq i \leq r-1$, $f(e_{0,2i}) = 1+i$ and $f(e_{0,2i+1}) = -(1+i)$; for $1 \leq i \leq s-1$, $f(e_{0,2i+1}) = r+i$ and $f(e_{0,2i+2}) = -(r+i)$; and for $1 \leq i \leq t$, $f(e_{0,2i-1}) = r+s-1+i$ and $f(e_{0,2i}) = -(r+s-1+i)$. This edge labeling induces a vertex labeling $f^+$ from $V(G)$ to $\{ \pm1, \ldots, \pm(r+s+t+1) \}$ as follows: $f^+(v_0) = r+s+t+1$; for $1 \leq i \leq 2r-1$, $f^+(v_i) = f(e_{0,i})$; $f^+(v_{2r}) = -(r+s+t+1)$ and $f^+(v_{2r+1}) = r+s+t$; for $1 \leq i \leq 2s$, $f^+(v_{2r,i}) = f(e_{2r,i})$; and for $1 \leq i \leq 2t$, $f^+(v_{2r+1,i}) = f(e_{2r+1,i})$. \qed
\end{proof}
\\[\baselineskip]
\indent The following lemma provides our first (and only) examples of caterpillars of diameter $4$ that are not super edge-graceful. 

\begin{lem}
\label{lem:nonSEGcaterpillars}
If $j,a_{j+2} \geq 1$ are odd integers with $j=1$ or $a_{j+2}=1$, then the rooted tree $RT(0^j, 1, a_{j+2})$ is not super edge-graceful.
\end{lem}  

\begin{proof}
Assume, to the contrary, that the rooted tree $G = \linebreak RT(0^j,1,a_{j+2})$ for odd integers $j=2r+1$ and $a_{j+2}=2t+1$, with $r,t \geq 0$, is super edge-graceful. Then there exists a bijection $f$ from $E(G)$ to $\{ 0, \pm1, \ldots, \pm(r+t+2) \}$ where either $f(e_{0,2r+2}) = 0$ or $f(e_{0,2r+3}) = 0$ such that the induced vertex labeling $f^+$ from $V(G)$ to $\{ \pm1, \pm2, \ldots, \pm(r+t+3) \}$ is a bijection. Such a vertex labeling must label one vertex by $r+t+3$ and another by $-(r+t+3)$. For any edge $e$, $-(r+t+2) \leq f(e) \leq r+t+2$, implying that $f^+(u)=r+t+3$ and $f^+(v)=-f^+(u)$ for some distinct vertices $u,v \in \{ v_0, v_{2r+2}, v_{2r+3} \}$. 

\textit{Case 1.} Assume $r \geq 0$ and $t=0$. If $f(e_{0,2r+2}) = 0$, then $f^+(v_{2r+2}) = f(e_{2r+2,1}) = f^+(v_{2r+2,1})$. Similarly, if $f(e_{0,2r+3})=0$, then $f^+(v_{2r+3}) = f(e_{2r+3,1}) = f^+(v_{2r+3,1})$. Thus, $RT(0^{2r+1},1,1)$ for integers $r \geq 0$ is not super edge-graceful.

\textit{Case 2.} Assume $r=0$ and $t \geq 1$. If $f(e_{0,2})=0$, then \linebreak $f^+(v_2) = f(e_{2,1}) = f^+(v_{2,1})$, which is impossible. Thus, we must have $f(e_{0,3})=0$ and $f$ must label the remaining $2t+4$ edges from $\{ \pm1, \pm2, \ldots, \pm(t+2) \}$. Likewise, the induced vertex labeling $f^+$ must be defined so that $f^+(u)=t+3$ and $f^+(v)=-f^+(u)$ for some distinct vertices $u,v \in \{ v_0,v_2,v_3 \}$. 
\begin{itemize}
	\item \textit{Subcase 1.} Let $f^+(v_0)$ be either $t+3$ or $-(t+3)$, and $f^+(v_2)=-f^+(v_0)$. Without loss of generality, assume $f^+(v_0) = f(e_{0,1}) + f(e_{0,2}) = t+3$. If $f(e_{0,2}) < 0$, then $f(e_{0,1}) > t+3$, a contradiction. Then $f(e_{0,2}) > 0$. Since $f^+(v_2) = f(e_{0,2}) + f(e_{2,1}) = -(t+3)$, we find that $f(e_{2,1}) < -(t+3)$, another contradiction.  	
	\item \textit{Subcase 2.} Let $f^+(v_0)$ be either $t+3$ or $-(t+3)$, and $f^+(v_3)=-f^+(v_0)$. Without loss of generality, assume $f^+(v_0) = f(e_{0,1}) + f(e_{0,2}) = t+3$. Any super edge-graceful labeling of $G=RT(0,1,2t+1)$ has $f(E(G)) - \{ 0 \} \subset f^+(V(G))$, so for some vertex $v$ in $G$, $f^+(v) = f(e_{0,2})$. However, the only remaining unlabeled vertices are pendant or the vertex $v_2$. No pendant vertex could be labeled $f(e_{0,2})$ because the incident edge is certainly not labeled $f(e_{0,2})$, which means $f^+(v_2) = f(e_{0,2}) + f(e_{2,1}) = f(e_{0,2})$. This implies $f(e_{2,1}) = 0$, a contradiction. 
	\item \textit{Subcase 3.} Let $f^+(v_2)$ be either $t+3$ or $-(t+3)$, and $f^+(v_3)=-f^+(v_2)$. Without loss of generality, assume $f^+(v_2) = f(e_{0,2}) + f(e_{2,1}) = t+3$. Similar to the previous subcase, for some vertex $v$, $f^+(v) = f(e_{0,2})$, and the only possible choice for this vertex is $v_0$. This means that $f^+(v_0) = f(e_{0,1}) + f(e_{0,2}) = f(e_{0,2})$, implying $f(e_{0,1}) = 0$, which is a contradiction.	 
\end{itemize}
Thus, $RT(0,1,2t+1)$ for integers $t \geq 1$ is not super edge-graceful. \qed
\end{proof}

\begin{lem}
\label{lem:a_{j+1}=1}
For all odd integers $j,a_{j+2} \geq 3$, the rooted tree \linebreak $RT(0^j, 1, a_{j+2})$ is super edge-graceful.
\end{lem}

\begin{proof}
Let $j=2r+1$ and $a_{2r+3}=2t+1$ for integers $r,t \geq 1$. Then $G=RT(0^{2r+1}, 1, 2t+1)$ has $2(r+t+2)+1$ edges and $2(r+t+3)$ vertices. An edge-labeling bijection $f$ maps $E(G)$ to $\{ 0, \pm1, \ldots, \pm(r+t+2) \}$ as follows: $f(e_{0,1})=-1$, $f(e_{0,2})=-2$, $f(e_{0,3})=3$, $f(e_{0,2r+2})=1$, $f(e_{0,2r+3})=0$, $f(e_{2r+2,1}) = r+t+2$, $f(e_{2r+3,1})=2$, $f(e_{2r+3,2})=-3$, and $f(e_{2r+3,3}) = -(r+t+2)$; for $2 \leq i \leq r$, $f(e_{0,2i}) = i+2$ and $f(e_{0,2i+1})= -(i+2)$; and for integers $2 \leq i \leq t$, $f(e_{2r+3,2i}) = r+1+i$ and $f(e_{2r+3,2i+1}) = -(r+1+i)$. This edge labeling induces a vertex labeling $f^+$ from $V(G)$ to $\{ \pm1, \ldots, \pm(r+t+3) \}$ as follows: $f^+(v_0)=1$; for integers $1 \leq i \leq 2r+1$, $f^+(v_i) = f(e_{0,i})$; $f^+(v_{2r+2}) = r+t+3$ and $f^+(v_{2r+3}) = -(r+t+3)$; $f^+(v_{2r+2,1})= r+t+2$; and for integers $1 \leq i \leq 2t+1$, $f^+(v_{2r+3,i}) = f(e_{2r+3,i})$. \qed
\end{proof}

\begin{lem}
\label{lem:o,o,o}
For all odd integers $j \geq 1$ and $a_{j+2} \geq a_{j+1} \geq 3$, the rooted tree $RT(0^j, a_{j+1}, a_{j+2})$ is super edge-graceful.
\end{lem}

\begin{proof}
Let $j=2r+1$, $a_{2r+2}=2s+1$, and $a_{2r+3}=2t+1$ for integers $r \geq 0$, $t \geq s \geq 1$. Then $G=RT(0^{2r+1}, 2s+1, 2t+1)$ has $2(r+s+t+2)+1$ edges and $2(r+s+t+3)$ vertices. An edge-labeling bijection $f$ maps $E(G)$ to $\{ 0, \pm1, \ldots, \pm(r+s+t+2) \}$ as follows: $f(e_{0,1}) = r+s+t+2$, $f(e_{0,2r+2}) = 1$, $f(e_{0,2r+3}) = 0$, $f(e_{2r+2,1}) = -1$, $f(e_{2r+2,2}) = -2$, $f(e_{2r+2,3}) = 3$, $f(e_{2r+3,1}) = 2$, $f(e_{2r+3,2}) = -3$, and $f(e_{2r+3,3}) = -(r+s+t+2)$; for $1 \leq i \leq r$, $f(e_{0,2i}) = 3+i$ and $f(e_{0,2i+1}) = -(3+i)$; for $2 \leq i \leq s$, $f(e_{2r+2,2i}) = r+2+i$ and $f(e_{2r+2,2i+1}) = -(r+2+i)$; and for $2 \leq i \leq t$, $f(e_{2r+3,2i}) = r+s+1+i$ and $f(e_{2r+3,2i+1}) = -(r+s+1+i)$. This edge labeling induces a vertex labeling $f^+$ from $V(G)$ to $\{ \pm1, \ldots, \pm(r+s+t+3) \}$ as follows: $f^+(v_0) = r+j+k+3$; for integers $1 \leq i \leq 2r+2$, $f^+(v_i) = f(e_{0,i})$; $f^+(v_{2r+3}) = -(r+s+t+3)$; for integers $1 \leq i \leq 2s+1$, $f^+(v_{2r+2,i}) = f(e_{2r+2,i})$; and for integers $1 \leq i \leq 2t+1$, $f^+(v_{2r+3,i}) = f(e_{2r+3,i})$. \qed
\end{proof}
\\[\baselineskip]
\indent Lemmas~\ref{lem:e,oppparity}, \ref{lem:o,e,e}, \ref{lem:nonSEGcaterpillars}, \ref{lem:a_{j+1}=1}, and~\ref{lem:o,o,o} imply the following theorem:

\begin{thm}
\label{thm:oddcaterpillars}
All caterpillars of diameter $4$ with odd size are super edge-graceful, except $RT(0^j,1,a_{j+2})$ where $j,a_{j+2}$ are positive odd integers with $j=1$ or $a_{j+2}=1$.
\end{thm}

\begin{ex}
Figure~\ref{fig:RT(0^3,2,4)-RT(0^3,3,5)} shows super edge-graceful labelings of two odd size caterpillars of diameter $4$, $RT(0^3,2,4)$ and $RT(0^3,3,5)$.
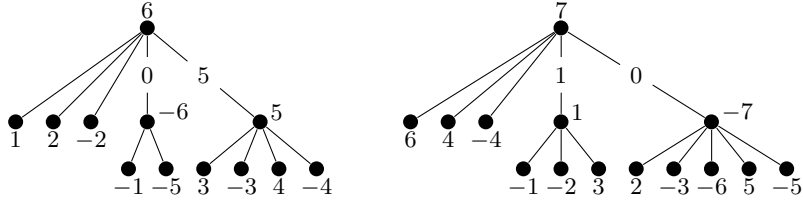
\begin{figure}[ht]
\centering
\begin{tikzpicture}[scale=2.5]
\tikzstyle{vertex}=[circle,fill=black,inner sep=2pt]

\node[vertex] (c-0) at (1,1) {};

\coordinate[label=above:\footnotesize $6$] (C-0) at (1,1);

\foreach \cent/\xcoord in {1/0.3, 2/0.5, 3/0.7, 4/1, 5/1.6}
	\node[vertex] (c-\cent) at (\xcoord,0.5) {};

\foreach \cent/\lab/\xcoord in {1/1/0.3, 2/2/0.5, 3/-2/0.7}
  \coordinate[label=below:\footnotesize $\lab$] (C-\cent) at (\xcoord,0.5);
  
\node[] (e-04) at (1,0.75) {\footnotesize $0$};
\node[] (e-05) at (1.3,0.75) {\footnotesize $5$};
  
\coordinate[label=right:\footnotesize $-6$] (C-4) at (1,0.55);
\coordinate[label=right:\footnotesize $5$] (C-5) at (1.6,0.55);

\foreach \pen/\lab/\xcoord in {41/-1/0.9, 42/-5/1.1, 51/3/1.3, 52/-3/1.5, 53/4/1.7, 54/-4/1.9}
{
 \node[vertex] (x-\pen) at (\xcoord,0.25) {};
 \coordinate[label=below:\footnotesize $\lab$] (X-\pen) at (\xcoord,0.25);
}

\foreach \cent in {1,2,3}
 \draw (c-0) -- (c-\cent);
\draw
(c-0) -- (e-04) -- (c-4)
(c-0) -- (e-05) -- (c-5);
\foreach \cent/\pen in {4/41, 4/42, 5/51, 5/52, 5/53, 5/54}
 \draw (c-\cent) -- (x-\pen); 

\node[vertex,xshift=6cm] (c-0) at (0.8,1) {};

\coordinate[label=above:\footnotesize $7$,xshift=6cm] (C-0) at (0.8,1);

\foreach \cent/\xcoord in {1/0, 2/0.2, 3/0.4, 4/0.8, 5/1.6}
	\node[vertex,xshift=6cm] (c-\cent) at (\xcoord,0.5) {};

\foreach \cent/\lab/\xcoord in {1/6/0, 2/4/0.2, 3/-4/0.4}
  \coordinate[label=below:\footnotesize $\lab$,xshift=6cm] (C-\cent) at (\xcoord,0.5);
  
\node[xshift=6cm] (e-04) at (0.8,0.75) {\footnotesize $1$};
\node[xshift=6cm] (e-05) at (1.2,0.75) {\footnotesize $0$};
  
\coordinate[label=right:\footnotesize $1$,xshift=6cm] (C-4) at (0.8,0.55);
\coordinate[label=right:\footnotesize $-7$,xshift=6cm] (C-5) at (1.6,0.55);

\foreach \pen/\lab/\xcoord in {41/-1/0.6, 42/-2/0.8, 43/3/1, 51/2/1.2, 52/-3/1.4, 53/-6/1.6, 54/5/1.8, 55/-5/2}
{
 \node[vertex,xshift=6cm] (x-\pen) at (\xcoord,0.25) {};
 \coordinate[label=below:\footnotesize $\lab$,xshift=6cm] (X-\pen) at (\xcoord,0.25);
}

\foreach \cent in {1,2,3}
 \draw (c-0) -- (c-\cent);
\draw
(c-0) -- (e-04) -- (c-4)
(c-0) -- (e-05) -- (c-5);
\foreach \cent/\pen in {4/41, 4/42, 4/43, 5/51, 5/52, 5/53, 5/54, 5/55}
 \draw (c-\cent) -- (x-\pen);

\end{tikzpicture}
\caption{$RT(0^3,2,4)$ and $RT(0^3,3,5)$ are super edge-graceful.}
\label{fig:RT(0^3,2,4)-RT(0^3,3,5)}
\end{figure}
\end{ex}

\section{Lobsters of diameter four}

\begin{figure}[ht]
\centering
\begin{tikzpicture}[scale=3.0]
\tikzstyle{vertex}=[circle,fill=black,inner sep=2pt]

\node[vertex] (c-0) at (0.4,0.75) {};

\foreach \cent/\xcoord in {1/0, 2/0.35, 3/0.8}
	\node[vertex] (c-\cent) at (\xcoord,0.5) {};

\foreach \pen/\xcoord in {11/-0.1, 12/0.1, 21/0.35, 31/0.6, 32/0.8, 33/1.0}
 	\node[vertex] (x-\pen) at (\xcoord,0.25) {};

\foreach \cent in {1,2,3}
 	\draw (c-0) -- (c-\cent);
\foreach \cent/\pen in {1/11, 1/12, 2/21, 3/31, 3/32, 3/33}
 	\draw (c-\cent) -- (x-\pen);
 
 
\node[vertex,xshift=4.5cm] (c-0) at (0.85,0.75) {};

\foreach \cent/\xcoord in {1/0, 2/0.2, 3/0.45, 4/0.8, 5/1.2, 6/1.7}
	\node[vertex,xshift=4.5cm] (c-\cent) at (\xcoord,0.5) {};

\foreach \pen/\xcoord in {31/0.4, 32/0.5, 41/0.7, 42/0.8, 43/0.9, 51/1.1, 52/1.2, 53/1.3, 61/1.5, 62/1.6, 63/1.7, 64/1.8, 65/1.9}
 	\node[vertex,xshift=4.5cm] (x-\pen) at (\xcoord,0.25) {};

\foreach \cent in {1,2,3,4,5,6}
 	\draw (c-0) -- (c-\cent);
\foreach \cent/\pen in {3/31, 3/32, 4/41, 4/42, 4/43, 5/51, 5/52, 5/53, 6/61, 6/62, 6/63, 6/64, 6/65}
 	\draw (c-\cent) -- (x-\pen); 
\end{tikzpicture}
\caption{The lobsters $RT(2,1,3)$ and $RT(0^2,2,3^2,5)$.}
\label{fig:gen-lobsters}
\end{figure}
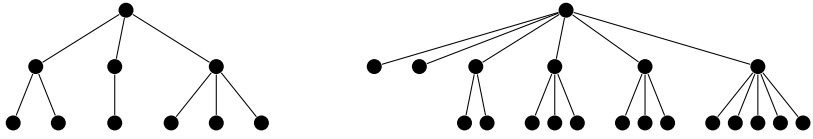

A lobster is a tree with the property that the removal of its endpoints leaves a caterpillar~\cite{JAGallian}. In this section, we examine rooted trees $RT(a_1, \ldots, a_n)$ that are lobsters of diameter $4$ and determine whether certain families of such trees are super edge-graceful. Recall that for integers $1 \leq i \leq n=j+k+l$, we have that $j$ is the cardinality of $\{ a_i : a_i=0 \}$, $k$ is the cardinality of $\{ a_i : a_i>0,~\mbox{even} \}$, and $l$ is the cardinality of $\{ a_i : a_i>0,~\mbox{odd} \}$. To ensure that a rooted tree is a lobster of diameter $4$, we assume in this section that $k+l \geq 3$. Figure~\ref{fig:gen-lobsters} shows the rooted trees $RT(2,1,3)$ and $RT(0^2,2,3^2,5)$, both of which are lobsters of diameter $4$.

Rooted trees that are lobsters of diameter $4$ will have $q = n + \sum_{i=j+1}^{n} a_i$ edges and $p = q+1$ vertices. We know that for $j+1 \leq i \leq j+k$, $a_i = 2b_i >0$ is even, where $b_i \geq 1$ is an integer, implying $\sum_{i=j+1}^{j+k} a_i = 2 \sum_{i=j+1}^{j+k} b_i$. Similarly, for $j+k+1 \leq i \leq n$, we have $a_i = 2b_i + 1 > 0$ is odd, where $b_i \geq 0$ is an integer, so that $\sum_{i=j+k+1}^{n} a_i = 2 \sum_{i=j+k+1}^{n} b_i + l$. Thus, the number of edges in a rooted tree that is a lobster is $q = j + k + 2l + 2\sum_{i=j+1}^{n} b_i$ where $b_i$ is an integer with $b_i \geq 1$ for $j+1 \leq i \leq j+k$ and $b_i \geq 0$ for $j+k+1 \leq i \leq n$. As in the previous section concerning caterpillars, we examine the rooted trees that are lobsters based on the parity of $q$, the size of the lobster.

\subsection{Even size lobsters of diameter four}

When $n = j+k+l$ and $l$ are of the same parity, that is, when $j,k$ are of the same parity, $q$ is even and $p$ is odd. 

\begin{lem}
\label{lem:j,k,l-odd}
If $n = j+k+l$ where $j,k,l \geq 1$ are odd integers with $k+l \geq 3$, then the rooted tree $RT(a_1, \ldots, a_n)$ is super edge-graceful.
\end{lem}

\begin{proof}
Let $j=2r-1$, $k=2s+1$, and $l=2t+1$ for integers $r \geq 1$ and $s,t \geq 0$ with $s+t \geq 1$ and let $G$ be the rooted tree defined in the statement of the lemma. Then $G$ has $2(r+s + 2t+1 + \sum_{i=2r}^{n} b_i)$ edges. An edge-labeling bijection $f$ maps $E(G)$ to $\{ \pm1, \ldots, \pm(r+s + 2t+1 + \sum_{i=2r}^{n} b_i) \}$ as follows: $f(e_{0,2(r+s)+1}) = 1$; for $1 \leq i \leq t$, $f(e_{0,2(r+s+i)}) = -(2i-1)$ and $f(e_{0,2(r+s+i)+1}) = 2i+1$; $f(e_{2(r+s+t)+1,1}) = -(2t+1)$; for $1 \leq i \leq t$, $f(e_{2(r+s+i)-1,1}) = -2(t+1-i)$ and $f(e_{2(r+s+i),1}) = 2(t+1-i)$; for $1 \leq i \leq r+s$, $f(e_{0,2i-1}) = 2t+1+i$ and $f(e_{0,2i}) = -(2t+1+i)$; for $2r \leq i \leq 2(r+s)$ and $1 \leq m \leq b_i$, $f(e_{i,2m-1}) = r+s+2t+1+m+\sum_{z=2r}^{i-1} b_z$ and $f(e_{i,2m}) = -f(e_{i,2m-1})$; for $2(r+s)+1 \leq i \leq n$ and $1 \leq m \leq b_i$, $f(e_{i,2m}) = r+s+2t+1+m+\sum_{z=2r}^{i-1} b_z$ and $f(e_{i,2m+1}) = -f(e_{i,2m})$.  \qed
\end{proof} 

\begin{lem}
\label{lem:j,k-even,l-odd}
If $n = j+k+l$ where $j,k \geq 0$ are even integers and $l \geq 1$ is an odd integer with $k+l \geq 3$, then the rooted tree $RT(a_1, \ldots, a_n)$ is super edge-graceful.
\end{lem}

\begin{proof}
Let $j=2r$, $k=2s$, and $l=2t+1$ for integers $r,s,t \geq 0$ with $s+t \geq 1$ and let $G$ be the rooted tree defined in the statement of the lemma. Then $G$ has $2(r+s + 2t+1 + \sum_{i=2r+1}^{n} b_i)$ edges. An edge-labeling bijection $f$ maps $E(G)$ to $\{ \pm1, \ldots, \pm(r+s + 2t+1 + \sum_{i=2r+1}^{n} b_i) \}$ as follows: $f(e_{0,2(r+s)+1}) = 1$; for $1 \leq i \leq t$, $f(e_{0,2(r+s+i)}) = -(2i-1)$ and $f(e_{0,2(r+s+i)+1}) = 2i+1$; $f(e_{2(r+s+t)+1,1}) = -(2t+1)$; for $1 \leq i \leq t$, $f(e_{2(r+s+i)-1,1}) = -2(t+1-i)$ and $f(e_{2(r+s+i),1}) = 2(t+1-i)$; for $1 \leq i \leq r+s$, $f(e_{0,2i-1}) = 2t+1+i$ and $f(e_{0,2i}) = -(2t+1+i)$; for $2r+1 \leq i \leq 2(r+s)$ and $1 \leq m \leq b_i$, $f(e_{i,2m-1}) = r+s+2t+1+m+\sum_{z=2r+1}^{i-1} b_z$ and $f(e_{i,2m}) = -f(e_{i,2m-1})$; for $2(r+s)+1 \leq i \leq n$ and $1 \leq m \leq b_i$, $f(e_{i,2m}) = r+s+2t+1+m+\sum_{z=2r+1}^{i-1} b_z$ and $f(e_{i,2m+1}) = -f(e_{i,2m})$.  \qed
\end{proof}

\begin{lem}
\label{lem:j,k,l-even}
If $n = j+k+l$ where $j,k,l \geq 0$ are even integers with $k+l \geq 3$, then the rooted tree $RT(a_1, \ldots, a_n)$ is super edge-graceful.
\end{lem}

\begin{proof}
Let $j=2r$, $k=2s$, and $l=2t$ for integers $r,s,t \geq 0$ with $s+t \geq 2$ and let $G$ be the rooted tree defined in the statement of the lemma. Then $G$ has $2(r+s + 2t + \sum_{i=2r+1}^{n} b_i)$ edges. An edge-labeling bijection $f$ maps $E(G)$ to $\{ \pm1, \ldots, \pm(r+s + 2t + \sum_{i=2r+1}^{n} b_i) \}$ as follows: for $1 \leq i \leq t$, $f(e_{0,2(r+s+i)-1}) = 2i-1$, $f(e_{0,2(r+s+i)}) = -(2i-1)$, $f(e_{2(r+s+i)-1,1}) = -2(t+1-i)$, and $f(e_{2(r+s+i),1}) = 2(t+1-i)$; for $1 \leq i \leq r+s$, $f(e_{0,2i-1}) = 2t+i$ and $f(e_{0,2i}) = -(2t+i)$; for $2r+1 \leq i \leq 2(r+s)$ and $1 \leq m \leq b_i$, $f(e_{i,2m-1}) = r+s+2t+m+\sum_{z=2r+1}^{i-1} b_z$ and $f(e_{i,2m}) = -f(e_{i,2m-1})$; for $2(r+s)+1 \leq i \leq n$ and $1 \leq m \leq b_i$, $f(e_{i,2m}) = r+s+2t+m+\sum_{z=2r+1}^{i-1} b_z$ and $f(e_{i,2m+1}) = -f(e_{i,2m})$.  \qed
\end{proof}

\begin{lem}
\label{lem:j,k-odd,l-even}
If $n = j+k+l$ where $j,k \geq 1$ are odd integers and $l \geq 0$ is an even integer with $k+l \geq 3$, then the rooted tree $RT(a_1, \ldots, a_n)$ is super edge-graceful.
\end{lem}

\begin{proof}
Let $j=2r-1$, $k=2s+1$, $l=2t$ for integers $r \geq 1$ and $s,t \geq 0$ with $s+t \geq 1$ and let $G$ be the rooted tree defined in the statement of the lemma. Then $G$ has $2(r+s + 2t + \sum_{i=2r}^{n} b_i)$ edges. An edge-labeling bijection $f$ maps $E(G)$ to $\{ \pm1, \ldots, \pm(r+s + 2t + \sum_{i=2r}^{n} b_i) \}$ as follows: for $1 \leq i \leq t$, $f(e_{0,2(r+s+i)-1}) = 2i-1$, $f(e_{0,2(r+s+i)}) = -(2i-1)$, $f(e_{2(r+s+i)-1,1}) = -2(t+1-i)$, and $f(e_{2(r+s+i),1}) = 2(t+1-i)$; for $1 \leq i \leq r+s$, $f(e_{0,2i-1}) = 2t+i$ and $f(e_{0,2i}) = -(2t+i)$; for $2r \leq i \leq 2(r+s)$ and $1 \leq m \leq b_i$, $f(e_{i,2m-1}) = r+s+2t+m+\sum_{z=2r}^{i-1} b_z$ and $f(e_{i,2m}) = -f(e_{i,2m-1})$; for $2(r+s)+1 \leq i \leq n$ and $1 \leq m \leq b_i$, $f(e_{i,2m}) = r+s+2t+m+\sum_{z=2r}^{i-1} b_z$ and $f(e_{i,2m+1}) = -f(e_{i,2m})$.  \qed
\end{proof}
\\[\baselineskip]
\indent Lemmas~\ref{lem:j,k,l-odd}, \ref{lem:j,k-even,l-odd}, \ref{lem:j,k,l-even}, and~\ref{lem:j,k-odd,l-even} imply the following theorem:

\begin{thm}
\label{thm:evenlobsters}
All lobsters of diameter $4$ with even size are super edge-graceful.
\end{thm}

\begin{ex}
Figure~\ref{fig:RT(0,2,3^2,5)} shows an even size lobster of diameter $4$, $RT(0,2,3^2,5)$, with a super edge-graceful labeling. 
\begin{figure}[ht]
\centering
\begin{tikzpicture}[scale=3.0]
\tikzstyle{vertex}=[circle,fill=black,inner sep=2pt]

\node[vertex] (c-0) at (1.2,1) {};

\coordinate[label=above:\footnotesize $3$] (C-0) at (1.2,1);

\node (e-02) at (0.8,0.75) {\footnotesize $-4$};
\node (e-03) at (1.1,0.75) {\footnotesize $1$};
\node (e-04) at (1.45,0.75) {\footnotesize $-1$};
\node (e-05) at (1.9,0.75) {\footnotesize $3$};

\foreach \cent/\xcoord in {1/-0.1, 2/0.4, 3/1, 4/1.7, 5/2.6}
	\node[vertex] (c-\cent) at (\xcoord,0.5) {};

\coordinate[label=below:\footnotesize $4$] (C-1) at (-0.1,0.5);
\coordinate[label=left:\footnotesize $-4$] (C-2) at (0.4,0.5);
\coordinate[label=left:\footnotesize $-1$] (C-3) at (1,0.5);
\coordinate[label=right:\footnotesize $1$] (C-4) at (1.7,0.5);
\coordinate[label=right:\footnotesize $0$] (C-5) at (2.6,0.525);

\foreach \pen/\lab/\xcoord in {21/5/0.3, 22/-5/0.5, 31/-2/0.8, 32/6/1, 33/-6/1.2, 41/2/1.5, 42/7/1.7, 43/-7/1.9, 51/-3/2.2, 52/8/2.4, 53/-8/2.6, 54/9/2.8, 55/-9/3}
{
 	\node[vertex] (x-\pen) at (\xcoord,0.25) {};
 	\coordinate[label=below:\footnotesize $\lab$] (C-1) at (\xcoord,0.25);
}

\draw (c-0) -- (c-1);

\foreach \edge/\cent in {02/2,03/3,04/4,05/5}
 	\draw (c-0) -- (e-\edge) -- (c-\cent);
\foreach \cent/\pen in {2/21, 2/22, 3/31, 3/32, 3/33, 4/41, 4/42, 4/43, 5/51, 5/52, 5/53, 5/54, 5/55}
 	\draw (c-\cent) -- (x-\pen); 
\end{tikzpicture}
\caption{$RT(0,2,3^2,5)$ is super edge-graceful.}
\label{fig:RT(0,2,3^2,5)}
\end{figure}
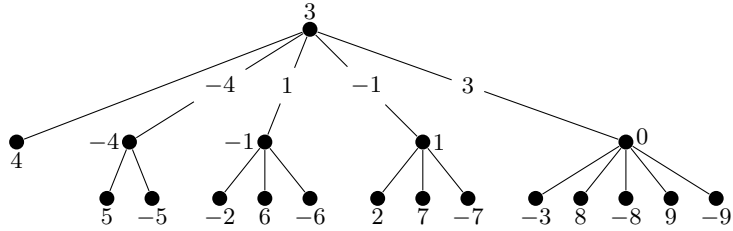  
\end{ex}

\subsection{Odd size lobsters of diameter four}

When $n = j+k+l$ and $l$ are of opposite parity, that is, when $j,k$ are of opposite parity, $q$ is odd and $p$ is even. If the odd size lobster $G=RT(a_1,\ldots,a_n)$ of diameter $4$ is to be super edge-graceful, then one of the edges in $E(G)$ must be labeled by $0$. Edges in $E(G) - \{ e_{0,j+1}, \ldots, e_{0,n} \}$ are incident to pendant vertices, so none of these edges can be labeled by $0$, implying $f(e_{0,m})=0$ for some \linebreak $m \in \{ j+1, \ldots, n \}$. 

\begin{lem}
\label{lem:j-even,k,l-odd}
If $n = j+k+l$ where $j \geq 0$ is an even integer and $k,l \geq 1$ are odd integers with $k+l \geq 3$, then the rooted tree $RT(a_1, \ldots, a_n)$ is super edge-graceful.
\end{lem}

\begin{proof}
Let $j=2r$, $k=2s-1$, $l=2t+1$ for integers $r,t \geq 0$ and $s \geq 1$ with $s+t \geq 2$ and let $G$ be the rooted tree defined in the statement of the lemma. Then $G$ has $q = 2(r+s + 2t + \sum_{i=2r+1}^{n} b_i)+1$ edges. 

\textit{Case 1.} If $l=1$ ($t=0$), then $k \geq 3$ ($s \geq 2$) and an edge-labeling bijection $f$ maps $E(G)$ to $\{ 0, \pm1, \ldots, \pm(r+s + \sum_{i=2r+1}^{n} b_i) \}$ as follows: $f(e_{0,2r+1}) = 0$, $f(e_{0,2(r+s)}) = 1$, $f(e_{2r+1,1}) = -1$, $f(e_{2(r+s),1}) = r+s + \sum_{i=2r+1}^{n} b_i$, and $f(e_{2r+1,2}) = -f(e_{2(r+s),1})$; for $1 \leq i \leq r$, $f(e_{0,2i-1}) = 1+i$ and $f(e_{0,2i}) = -(1+i)$; for $1 \leq i \leq s-1$, $f(e_{0,2(r+i)}) = r+1+i$ and $f(e_{0,2(r+i)+1}) = -(r+1+i)$; for $2 \leq m \leq b_{2r+1}$, $f(e_{2r+1,2m-1}) = r+s-1+m$ and $f(e_{2r+1,2m}) = -f(e_{2r+1,2m-1})$; for $2r+2 \leq i \leq 2(r+s)-1$ and $1 \leq m \leq b_i$, $f(e_{i,2m-1}) = r+s-1+m+\sum_{z=2r+1}^{i-1} b_z$ and $f(e_{i,2m}) = -f(e_{i,2m-1})$; for $2(r+s) \leq i \leq n$ and $1 \leq m \leq b_i$, $f(e_{i,2m}) = r+s-1+m+\sum_{z=2r+1}^{i-1} b_z$ and $f(e_{i,2m+1}) = -f(e_{i,2m})$.

\textit{Case 2.} If $l \geq 3$ ($t \geq 1$), then $k \geq 1$ ($s \geq 1$) and an edge-labeling bijection $f$ maps $E(G)$ to $\{ 0, \pm1, \ldots, \pm(r+s + 2t + \sum_{i=2r+1}^{n} b_i) \}$ as follows: $f(e_{0,2r+1}) = 0$; for $1 \leq i \leq t$, $f(e_{0,2(r+s+i-1)}) = 2i-1$ and $f(e_{0,2(r+s+i)-1}) = -(2i-1)$; $f(e_{0,2(r+s+t)}) = 2t+1$, $f(e_{2r+1,1}) = -(2t+1)$, $f(e_{2(r+s+t),1}) = -2$, $f(e_{2r+1,2}) = 2$, $f(e_{2(r+s),1}) = r+s + 2t + \sum_{i=2r+1}^{n} b_i$, $f(e_{2(r+s)+1,1}) = -f(e_{2(r+s),1})$; for $1 \leq i \leq t-1$, $f(e_{2(r+s+i),1}) = -2(t-i+1)$ and $f(e_{2(r+s+i)+1,1}) = 2(t-i+1)$; for $1 \leq i \leq r$, $f(e_{0,2i-1}) = 2t+1+i$ and $f(e_{0,2i}) = -(2t+1+i)$; for $1 \leq i \leq s-1$, $f(e_{0,2(r+i)}) = r+2t+1+i$ and $f(e_{0,2(r+i)+1}) = -f(e_{0,2(r+i)})$; for $2 \leq m \leq b_{2r+1}$, $f(e_{2r+1,2m-1}) = r+s+2t-1+m$ and $f(e_{2r+1,2m}) = -f(e_{2r+1,2m-1})$; for $2r+2 \leq i \leq 2(r+s)-1$ and $1 \leq m \leq b_i$, $f(e_{i,2m-1}) = r+s+2t-1+m+\sum_{z=2r+1}^{i-1} b_z$ and $f(e_{i,2m}) = -f(e_{i,2m-1})$; for $2(r+s) \leq i \leq n$ and $1 \leq m \leq b_i$, $f(e_{i,2m}) = r+s+2t-1+m+\sum_{z=2r+1}^{i-1} b_z$ and $f(e_{i,2m+1}) = -f(e_{i,2m})$. \qed
\end{proof}

\begin{lem}
\label{lem:j,l-even,k-odd}
If $n = j+k+l$ where $j,l \geq 0$ are even integers and $k \geq 3$ is an odd integer, then the rooted tree $RT(a_1, \ldots, a_n)$ is super edge-graceful.
\end{lem}

\begin{proof}
Let $j=2r$, $k=2s+1$, $l=2t$ for integers $r,t \geq 0$ and $s \geq 1$, and let $G$ be the rooted tree defined in the statement of the lemma. Then $G$ has $q = 2(r+s + 2t + \sum_{i=2r+1}^{n} b_i)+1$ edges. 

\textit{Case 1.} If $l=0$ ($t=0$), then an edge-labeling bijection $f$ maps $E(G)$ to $\{ 0, \pm1, \ldots, \pm(r+s + \sum_{i=2r+1}^{n} b_i) \}$ as follows: $f(e_{0,2r+1}) = 0$, $f(e_{0,2r+2}) = 1$, $f(e_{0,2r+3}) = r+s + \sum_{i=2r+1}^{n} b_i$, $f(e_{2r+1,1}) = -1$, and $f(e_{2r+1,2}) =  -f(e_{0,2r+3})$; for $1 \leq i \leq r$, $f(e_{0,2i-1}) = 1+i$ and $f(e_{0,2i}) = -(1+i)$; for $2 \leq i \leq s$, $f(e_{0,2(r+i)}) = r+i$ and $f(e_{0,2(r+i)+1}) = -(r+i)$; for $2 \leq m \leq b_{2r+1}$, $f(e_{2r+1,2m-1}) = r+s-1+m$ and $f(e_{2r+1,2m}) = -(r+s-1+m)$; for $2r+2 \leq i \leq 2(r+s)+1$ and $1 \leq m \leq b_i$, $f(e_{i,2m-1}) = r+s-1+m+\sum_{z=2r+1}^{i-1} b_z$ and $f(e_{i,2m}) = -f(e_{i,2m-1})$.

\textit{Case 2.} If $l \geq 2$ ($t \geq 1$), then $f$ maps $E(G)$ to $\{ 0, \pm1, \ldots, \pm(r+s + 2t + \sum_{i=2r+1}^{n} b_i) \}$ as follows: $f(e_{0,2r+1})=0$, $f(e_{0,2r+2}) = 2$, $f(e_{0,2r+3}) = -(2t+1)$, $f(e_{2r+1,1}) = -2$, $f(e_{2r+1,2}) = 2t+1$, \linebreak $f(e_{2(r+s+1),1})=r+s + 2t + \sum_{i=2r+1}^{n} b_i$, and $f(e_{2(r+s+1)+1,1}) = -f(e_{2(r+s+1),1})$; for $1 \leq i \leq r$, $f(e_{0,2i-1}) = 2t+1+i$ and $f(e_{0,2i})=-(2t+1+i)$; for $2 \leq i \leq s$, $f(e_{0,2(r+i)})=2t+r+i$ and $f(e_{0,2(r+i)+1}) = -(2t+1+i)$; for $1 \leq i \leq t$, $f(e_{0,2(r+s+i)})=2i-1$, $f(e_{0,2(r+s+i)+1})=-(2i-1)$; for $2 \leq i \leq t$, $f(e_{2(r+s+i),1})=-2(t+2-i)$, \linebreak $f(e_{2(r+s+i)+1,1})=2(t+2-i)$; for $2 \leq m \leq b_{2r+1}$, $f(e_{2r+1,2m-1}) = 2t+r+s-1+m$ and $f(e_{2r+1,2m}) = -f(e_{2r+1,2m-1})$; for $2r+2 \leq i \leq 2(r+s)+1$ and $1 \leq m \leq b_i$, $f(e_{i,2m-1}) = 2t+r+s-1+m+\sum_{z=2r+1}^{i-1} b_z$ and $f(e_{i,2m}) = -f(e_{i,2m-1})$; for $2(r+s+1) \leq i \leq n$ and $1 \leq m \leq b_i$, $f(e_{i,2m}) = 2t+r+s-1+m+\sum_{z=2r+1}^{i-1} b_z$ and $f(e_{i,2m+1}) = -f(e_{i,2m})$. \qed
\end{proof}

\begin{lem}
\label{lem:badlobsters}
For all odd integers $j \geq 1$ and integers $l \geq 3$, the rooted tree $RT(0^j,1^l)$ is not super edge-graceful.
\end{lem}

\begin{proof}
Assume, to the contrary, that the rooted tree $RT(0^j,1^l)$ for odd integers $j \geq 1$ and any integer $l \geq 3$ is super edge-graceful. Such a rooted tree has odd size so some edge must be labeled by $0$; however, no edge incident with a pendant vertex may be labeled by $0$. Then $f(e_{0,m})=0$ for some integer $m \in \{ j+1, \ldots, j+l \}$. This implies $f^+(v_m)=f(e_{m,1})=f^+(v_{m,1})$, contradicting the definition of a super edge-graceful labeling. \qed
\end{proof}

\begin{lem}
\label{lem:j-odd,k,l-even}
If $n = j+k+l$ where $j \geq 1$ is an odd integer, $k \geq 2$ is an even integer, and $l=1$ or $l=2$, then the rooted tree $RT(a_1, \ldots, a_n)$ is super edge-graceful.
\end{lem}

\begin{proof}
Let $j=2r+1$ for integers $r \geq 0$ and $k=2s$ for integers $s \geq 1$. 

\textit{Case 1.} Suppose $l=1$ and let $G$ be the rooted tree defined in the statement of the lemma. Then $G$ has $q = 2(r+s +1 + \sum_{i=2r+2}^{n} b_i)+1$ edges and an edge-labeling bijection $f$ maps $E(G)$ to $\{ 0, \pm1, \ldots, \pm(r+s + 1 + \sum_{i=2r+2}^{n} b_i) \}$ as follows: $f(e_{0,2r+2}) = 0$, $f(e_{0,n}) = 1$, $f(e_{2r+2,1}) = -1$, $f(e_{n,1}) = r+s + 1 + \sum_{i=2r+2}^{n} b_i$, and $f(e_{2r+2,2}) = -f(e_{n,1})$; for $1 \leq i \leq r$, $f(e_{0,2i-1}) = 1+i$ and $f(e_{0,2i}) = -(1+i)$; $f(e_{0,2r+1}) = r+2$ and $f(e_{0,2r+3}) = -(r+2)$; for $2 \leq i \leq s$, $f(e_{0,2(r+i)}) = r+1+i$ and $f(e_{0,2(r+i)+1}) = -(r+1+i)$; for $2 \leq m \leq b_{2r+2}$, $f(e_{2r+2,2m-1}) = r+s+m$ and $f(e_{2r+2,2m}) = -(r+s+m)$; for $2r+3 \leq i \leq 2(r+s)+1$ and $1 \leq m \leq b_i$, $f(e_{i,2m-1}) = r+s+m+\sum_{z=2r+2}^{i-1} b_z$ and $f(e_{i,2m}) = -f(e_{i,2m-1})$; for $1 \leq m \leq b_n$, $f(e_{n,2m}) = r+s+m+\sum_{z=2r+2}^{2(r+s)+1} b_z$ and $f(e_{n,2m+1}) = -f(e_{n,2m})$.

\textit{Case 2.} Suppose $l=2$ and let $G$ be the rooted tree defined in the statement of the lemma. Then $G$ has $q = 2(r+s + 2 + \sum_{i=2r+2}^{n} b_i)+1$ edges and an edge-labeling bijection $f$ maps $E(G)$ to $\{ 0, \pm1, \ldots, \pm(r+s+ 2 + \sum_{i=2r+2}^{n} b_i) \}$ as follows: $f(e_{0,2r+2})=0$, $f(e_{0,1}) = 1$, $f(e_{2r+2,1}) = -1$, $f(e_{0,2(r+s+1)}) = 2$, $f(e_{0,2(r+s+1)+1}) = -2$, $f(e_{2r+3,1}) = 3$, $f(e_{2r+3,2}) = -3$, $f(e_{2(r+s+1),1}) = -4$, and $f(e_{2(r+s+1)+1,1}) = 4$, $f(e_{0,2r+3}) = r+s+ 2 + \sum_{i=2r+2}^{n} b_i$, and \linebreak $f(e_{2r+2,2}) = -f(e_{0,2r+3})$; for $1 \leq i \leq r$, $f(e_{0,2i}) = 4+i$ and $f(e_{0,2i+1}) = -(4+i)$; for $2 \leq i \leq s$, $f(e_{0,2i-1}) = r+3+i$, $f(e_{0,2i}) = -(r+3+i)$; for $2 \leq m \leq b_{2r+2}$, $f(e_{2r+2,2m-1}) = r+s+2+m$ and $f(e_{2r+2,2m}) = -f(e_{2r+2,2m-1})$; for $2 \leq m \leq b_{2r+3}$, $f(e_{2r+3,2m-1}) = r+s+b_{2r+2}+1+m$ and $f(e_{2r+3,2m}) = -f(e_{2r+3,2m-1})$; for $2r+4 \leq i \leq 2(r+s)+1$ and $1 \leq m \leq b_i$, $f(e_{i,2m-1}) = r+s+1+m+ \sum_{z=2r+2}^{i-1} b_z$ and $f(e_{i,2m}) = -f(e_{i,2m-1})$; for $i \in \{ 2r+2s+2, 2r+2s+3 \}$ and $1 \leq m \leq b_i$, $f(e_{i,2m}) = r+s+1+m+ \sum_{z=2r+2}^{i-1} b_z$ and $f(e_{i,2m+1}) = -f(e_{i,2m})$. \qed
\end{proof}

\begin{ex}
Figure~\ref{fig:RT(2,3^2,5)} shows an odd size lobster of diameter $4$, $RT(2,3^2,5)$, with a super edge-graceful labeling. 
\begin{figure}[ht]
\centering
\begin{tikzpicture}[scale=3.0]
\tikzstyle{vertex}=[circle,fill=black,inner sep=2pt]

\node[vertex] (c-0) at (1.2,1) {};

\coordinate[label=above:\footnotesize $3$] (C-0) at (1.2,1);

\node (e-02) at (0.8,0.75) {\footnotesize $0$};
\node (e-03) at (1.1,0.75) {\footnotesize $1$};
\node (e-04) at (1.45,0.75) {\footnotesize $-1$};
\node (e-05) at (1.9,0.75) {\footnotesize $3$};

\foreach \cent/\xcoord in {2/0.4, 3/1, 4/1.7, 5/2.6}
	\node[vertex] (c-\cent) at (\xcoord,0.5) {};

\coordinate[label=left:\footnotesize $-1$] (C-2) at (0.4,0.5);
\coordinate[label=left:\footnotesize $9$] (C-3) at (1,0.5);
\coordinate[label=right:\footnotesize $-9$] (C-4) at (1.7,0.5);
\coordinate[label=right:\footnotesize $1$] (C-5) at (2.6,0.525);

\foreach \pen/\lab/\xcoord in {21/-3/0.3, 22/2/0.5, 31/8/0.8, 32/4/1, 33/-4/1.2, 41/-8/1.5, 42/5/1.7, 43/-5/1.9, 51/-2/2.2, 52/6/2.4, 53/-6/2.6, 54/7/2.8, 55/-7/3}
{
 	\node[vertex] (x-\pen) at (\xcoord,0.25) {};
 	\coordinate[label=below:\footnotesize $\lab$] (C-1) at (\xcoord,0.25);
}

\foreach \edge/\cent in {02/2,03/3,04/4,05/5}
 	\draw (c-0) -- (e-\edge) -- (c-\cent);
\foreach \cent/\pen in {2/21, 2/22, 3/31, 3/32, 3/33, 4/41, 4/42, 4/43, 5/51, 5/52, 5/53, 5/54, 5/55}
 	\draw (c-\cent) -- (x-\pen); 
\end{tikzpicture}
\caption{$RT(2,3^2,5)$ is super edge-graceful.}
\label{fig:RT(2,3^2,5)}
\end{figure}
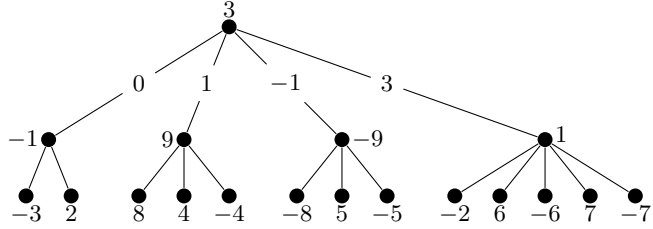  
\end{ex}

\subsubsection{Conjectures concerning odd size lobsters of diameter four}

We are unable to find general methods that describe super edge-graceful labelings for a few families of odd size lobsters of diameter $4$, although we are able to show that certain lobsters in these families are super-edge graceful.  

\begin{conj}
\label{conj:j,l-even,k-odd-bad}
If $n = j+k+l$ where $j \geq 0$ is and even integer, $k=1$, and $l \geq 2$ is an even integer, then the rooted tree $RT(a_1, \ldots, a_n)$ is super edge-graceful.
\end{conj}

\begin{conj}
\label{conj:j-odd,k,l-even-bad}
If $n = j+k+l$ where $j \geq 1$ is an odd integer, $k \geq 0$ is an even integer, and $l \geq 4$ is an even integer, then the rooted tree $RT(a_1, \ldots, a_n)$ is super edge-graceful provided $a_i \geq 3$ for some integer $i \in \{ j+k+1, \ldots , n \}$ when $k=0$.
\end{conj}

\begin{conj}
\label{conj:j,l-odd,k-even-bad}
If $n = j+k+l$ where $j \geq 1$ is an odd integer, $k \geq 0$ is an even integer, and $l \geq 3$ is an odd integer, then the rooted tree $RT(a_1, \ldots, a_n)$ is super edge-graceful provided $a_i \geq 3$ for some integer $i \in \{ j+k+1, \ldots , n \}$ when $k=0$.
\end{conj}

\bibliographystyle{plain}

\begin{thebibliography}{MMM}

\bibitem{CLGS} P.T.~Chung, S.M.~Lee, W.Y.~Gao, and K.~Schaffer, On the super edge graceful trees of even orders, \textit{Congr. Numer.}, \textbf{181} (2006) 5-17.

\bibitem{CFKX} S.~Cichacz, D.~Froncek, A.~Khodkar, and W.~Xu, Super edge-gradeful paths and cycles, \textit{Bull. Inst. Combin. Appl.}, \textbf{57} (2009) 79-90.

\bibitem{JAGallian} J.A.~Gallian, A dynamic survey of graph labeling, \textit{Elec. J. Combin.}, No. DS6, accessed November 19, 2010, {\tt http://www.combinatorics.org/Surveys/ds6.pdf}.

\bibitem{GWGolomb} S.W.~Golomb, How to number a graph, \textit{Graph Theory and Computing}, edited by R.C.~Read, Academic Press, N.Y. (1972) 23-37. 

\bibitem{HS} F.~Harary and A.~Schwenk, The number of caterpillars, \textit{Discrete Math.}, \textbf{6} (1973) 359-365.

\bibitem{LH} S.M.~Lee and Y.S.~Ho, All trees of odd order with three even vertices are super edge-graceful, \textit{J. Combin. Math. Combin. Comput.}, \textbf{62} (2007) 53-64.

\bibitem{LK} S.M.~Lee and M.C.~Kong, On super edge magic $n$-stars, \textit{J. Combin. Math. Combin. Comput.}, \textbf{42} (2002), 61-77.

\bibitem{SPLo} S.P.~Lo, On edge-graceful labelings of graphs, \textit{Congr. Numer.}, \textbf{50} (1985) 231-241.

\bibitem{MS} J.~Mitchem and A.~Simoson, On edge-graceful and super edge-graceful labelings of graphs, \textit{Ars Combin.}, \textbf{37} (1994) 97-111.

\bibitem{ARosa} A.~Rosa, On certain valuations of the vertices of a graph, \textit{Theory of Graphs (Internat. Symp., Rome, July 1966)}, Gordon and Breach, N.Y. and Dunod Paris (1967) 349-355.

\end{thebibliography}

\end{document}